\renewcommand{\baselinestretch}2

\documentclass[12pt]{amsart}
\usepackage{amsmath}
\usepackage{amsfonts}
\usepackage{amssymb}
\usepackage{graphicx, color}

{\theoremstyle{plain}%
  \newtheorem{theorem}{Theorem}[section]
  
\newtheorem{proposition}[theorem]{Proposition}
\newtheorem{lemma}[theorem]{Lemma}
}

{\theoremstyle{remark}

}

{\theoremstyle{definition}

\newtheorem{example}{Example}
}

\input xy
\xyoption{all}

\begin{document}
\title[Orthogonally additive holomorphic functions]{Orthogonally additive holomorphic functions of bounded type over $C(K)$}
\author{Daniel Carando}
\address{Depto. de Matem\'{a}tica, Pab. I, Fac. de Cs. Exactas y Naturales, Universidad
de Buenos Aires, (1428) Buenos Aires, Argentina}
\email{dcarando@dm.uba.ar}
\author{Silvia Lassalle}
\address{Depto de Matem\'{a}tica, Pab. I, Fac. de Cs. Exactas y Naturales, Universidad
de Buenos Aires, (1428) Buenos Aires, Argentina}
\email{slassall@dm.uba.ar}
\author{ Ignacio Zalduendo}
\address{Depto. de Matem\'{a}tica, Universidad Torcuato Di Tella,
Mi\~{n}ones 2177 (C1428ATG), Buenos Aires, Argentina}
\email{nacho@utdt.edu}

\thanks{Partially supported by PICT 05 33042. The first author 
was supported in part by UBACyT X038 and PICT 06-00587. The second author
was also supported by UBACyT X863 and PICT 06-00897. The three authors are members of CONICET}

\begin{abstract}
It is known that all $k$-homogeneous orthogonally additive polynomials $P$ over $C(K)$ are
of the form
$$
P(x)=\int_K x^k \, d\mu .
$$
Thus $x\mapsto x^k$ factors all orthogonally additive polynomials through some linear form $\mu$.
We show that no such linearization is possible without homogeneity. However, we also show that
every orthogonally additive holomorphic functions of bounded type $f$ over $C(K)$ is
of the form
$$
f(x)=\int_K h(x) \, d\mu
$$
for some $\mu$ and holomorphic $h\colon C(K) \rightarrow L^1(\mu)$ of bounded type.
\end{abstract}


\maketitle
\section*{Introduction.}

A holomorphic function $f\colon C(K)\rightarrow F$ is said to be orthogonally additive if
$f(x+y)=f(x)+f(y)$ whenever $x$ and $y$ are mutually orthognal (i.e., $xy=0$ over $K$).
In this paper we study orthogonal additivity of complex-valued bounded-type
holomorphic functions over $C(K)$.

Recall \cite{D3} that $f\colon E \rightarrow F$ is  of bounded type if it is bounded on all bounded
subsets of $E$. The Taylor series $f=\sum_k P_k$ of such a function have infinite radii of uniform convergence,
i.e.: $\limsup \Vert P_k \Vert^{\frac{1}{k}}=0$. We will denote the space of all such functions over $E$
by $\mathcal{H}_b (E,F)$, or simply $\mathcal{H}_b (E)$ if $F=\mathbb{C}$.

It was proved by \cite{PGV} and \cite{BLLl} (see also \cite{CaLaZa}) that each orthogonally
additive $k$-homogeneous polynomial $P$ over $C(K)$ is represented by a regular Borel measure
$\mu$ on $K$, in the following sense:
$$
P(x)=\int_K x^k \, d\mu \quad \text{ for all } x\in C(K).
$$
This is a linearization result; the `universal' polynomial
$h(x)=x^k$ linearizes all orthogonally additive $k$-homogeneous polynomials $P$,
$$
\xymatrix{
& C(K) \ar[d]_h \ar[r]^P & \mathbb{C} \\
& C(K) \ar[ur]_{\mu} & }
$$

The following lemma shows that $f$ is orthogonally additive if and only if all polynomials
in its Taylor series expansion at zero are orthogonally additive. Note that this is not the case
for expansions around other points.

\begin{lemma}
Let $f\colon C(K) \rightarrow \mathbb{C}$ be holomorphic, and $f=\sum_{k=0}^\infty P_k$ its Taylor series
at zero. Then $f$ is orthogonally additive if and only if all the $P_k$'s are orthogonally additive.
\end{lemma}

\begin{proof}
Say $x$ and $y$ are orthogonal. If $f$ is orthogonally additive, we may write
\begin{align*}
P_k(x+y)
& =\frac{1}{k!}\frac{d^k}{d\lambda}(0)\left(f(\lambda (x+y))\right) \\
& =\frac{1}{k!}\frac{d^k}{d\lambda}(0)\left(f(\lambda x)+ f(\lambda y)\right) \\
& =\frac{1}{k!}\frac{d^k}{d\lambda}(0) f(\lambda x) + \frac{1}{k!}\frac{d^k}{d\lambda}(0) f(\lambda y) \\
& =P_k(x) + P_k(y).
\end{align*}
For the converse, if all $P_k$'s are orthogonally additive,
\begin{align*}
f(x+y)&=\sum_{k=0}^\infty P_k(x+y)=\sum_{k=0}^\infty (P_k(x)+P_k(y))\\
      &=\sum_{k=0}^\infty P_k(x)+\sum_{k=0}^\infty P_k(y)\\
      &=f(x)+f(y).
\end{align*}
\end{proof}

Given the Lemma, one might expect that in the holomorphic setting
one will have a linearizing factorization where
the polynomial $h(x)=x^k$ is replaced by $h(x)=\Phi \circ x$, with
$\Phi \colon \mathbb{C} \rightarrow \mathbb{C}$ holomorphic and $\Phi(0)=0$
(for example $h(x)=e^x - 1$). The first objection to such a factorization
is that it would produce a linearization of the space of orthogonally additive
holomorphic functions through $C(K)$, and thus a Banach predual of the space; but spaces
of holomorphic functions are seldom Banach. A second objection to this line of thought
is that for any fixed $\Phi (z)=\sum_{k=1}^\infty a_k z^k$, the factorization can only
take place for functions $f$ with limited growth: indeed, if
$$
f(x)=\int_K \Phi \circ x \, d\mu = \int_K \sum_{k=1}^\infty a_k x^k \, d\mu
= \sum_{k=1}^\infty a_k \int_K x^k \, d\mu ,
$$
then
$$
|f(x)|\leq \sum_{k=1}^\infty |a_k| |\mu|(K) \Vert x \Vert^k = |\mu|(K) \sum_{k=1}^\infty |a_k| \Vert x \Vert^k ,
$$
so the growth of $f$ is bounded by that of $\Phi$.

We will show (Example 1) that there are stronger, algebraic,
obstructions and that such a factorization cannot be obtained even for non-homogeneous
polynomials of degree two.

In Section 1 we study the relationship between orthogonal additivity and integrality,
and prove that all orthogonally additive functions of bounded type are integral.
We also show an orthogonally additive polynomial of degree two that does not factor
through $C(K)$.

In Section 2 we characterize orthogonally additive holomorphic
functions of bounded type as those which may be written as
$$
f(x)=\int_K h(x) \, d\mu \quad \text{ for all } x \text{ in } C(K),
$$
for some measure $\mu$ on $K$ and $h\colon C(K) \rightarrow L^1(\mu)$, holomorphic and of bounded type.

\section{Orthogonal additivity and integrality}

Recall \cite{D1} that a $k$-homogeneous polynomial $P$ over a Banach space $E$ is
said to be integral if it can be represented by a Borel measure $\nu$ over
the unit ball $B_{E'}$ of the dual space (with the weak$^*$ topology) in the
following way
$$
P(x)=\int_{B_{E'}} \widehat{x}^k \, d\nu \quad \text{ for all } x\in E
$$
(where $\widehat{x}(\gamma)=\gamma(x)$). Recall also \cite{DiGaMaZa} that an integral holomorphic
function $f\colon B^\circ_E \rightarrow \mathbb{C}$ on the open unit ball
of $E$ is one which may be written as
$$
f(x)=\int_{B_{E'}} \frac{1}{1-\widehat{x}} \, d\nu \quad \text{ for all } x\in B^\circ_E\textcolor{red}{.}
$$
and some Borel measure $\nu$ on $B_{E'}$.

Note that when $E=C(K)$, the dual $E'$ is the space of regular Borel measures on $K$,
and that $K$ may be identified with part of the unit ball $B_{E'}$: the point measures
$\delta_a$, $a\in K$.

We begin with some observations on the case of orthogonally additive $k$-homogeneous polynomials.
Consider the closed subspace $X_k$ of $C(B_{E'})$ spanned by $\{\widehat{x}^k : x\in E\}$.
The space $X_k$ is then the symmetric $\epsilon$-tensor product $\widehat{\bigotimes}_{\epsilon ,k,s} E$
and integral $k$-homogeneous polynomials are those which linearize through $X_k$ \cite{D3}.
$$
\xymatrix{
x \ar@{|->}[d] & E \ar[d] \ar[r]^P & \mathbb{C} \\
{\widehat{x}^k } & X_k \ar[ur]^{\nu} &}
$$
Now if $E=C(K)$, we have the linear map $R_k\colon X_k\rightarrow C(K)$ given by
$R_k(\widehat{x}^k )=x^k$. It is easy to see that, since $K\subset B_{E'}$,
this map is well-defined and has norm one. The diagram
$$
\xymatrix{
x \ar@{|->}[d] & C(K) \ar[d] \ar[r]^P & \mathbb{C} \\
{\widehat{x}^k } & X_k \ar[d]_{R_k} \ar[ur]^{\nu} & \\
         & C(K) \ar[uur]_{\mu} &}
$$
then shows (applying \cite{PGV}, \cite{BLLl}) that all orthogonally additive $k$-homogeneous
polynomials are integral.
Indeed, if $P(x)=\mu(x^k)$, it is also $\nu(\widehat{x}^k)$.
This shows the well known fact that the `integrating' polynomial $P(x)=\int_K x(t)^k \, d\mu(t)$ on $C(K)$ is
integral.

It is also true that orthogonally additive holomorphic functions of bounded type are integral.
For the proof, we will need the following. Recall \cite{AB} that any $f\in \mathcal{H}_b(E)$
has a canonical extension to the bidual $E''$. This extension, $\overline{f}$, called the Aron-Berner extension
is also holomorphic and of bounded type.

Given a Borel subset $A$ of $K$, and $f\colon C(K) \rightarrow \mathbb{C}$ a bounded-type holomorphic
function, we define
$$
f_A(x)=\overline{f}(1_A x),
$$
where $\overline{f}$ is the Aron-Berner extension of $f$.

\begin{proposition}\label{bOA implies integral}
If $f\in \mathcal{H}_b (C(K))$ is orthogonally additive, then it is integral.
\end{proposition}

\begin{proof}
Let $f=\sum_{k=1}^\infty P_k$ be the Taylor series expansion of $f$ at zero, and note that each $P_k$
is orthogonally additive and hence integral. By \cite{DiGaMaZa} we need only verify that
$\sum_{k=1}^\infty \Vert P_k \Vert_I < \infty$. Since $f$ is of bounded type, we have
$\limsup \Vert P_k \Vert^{\frac{1}{k}}=0$, so $\sum_{k=1}^\infty \Vert P_k \Vert < \infty$.
But if $P$ is an orthogonally additive $k$-homogeneous polynomial,
then its natural and integral norms are equal:
In general, one has $\Vert P \Vert \leq \Vert P \Vert_I$; to see $\Vert P \Vert_I \leq \Vert P \Vert$,
given $\varepsilon > 0$, we find a measure $\nu$ on $B_{E'}$ representing $P$ and such that
$|\nu| < \Vert P \Vert + \varepsilon$: since $P$ is orthogonally additive, consider a measure $\mu$
on $K$ such that $P(x)=\int_K x^k \, d\mu$, and take $\nu=\mu\circ R_k$:
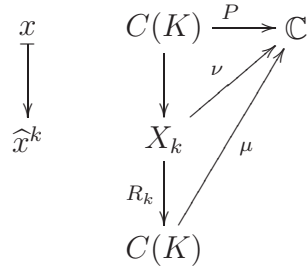
$$
\xymatrix{
x \ar@{|->}[d] & C(K) \ar[d] \ar[r]^P & \mathbb{C} \\
{\widehat{x}^k } & X_k \ar[d]_{R_k} \ar[ur]^{\nu} & \\
         & C(K) \ar[uur]_{\mu} &}
$$
Since $\Vert R_k \Vert=1$,
$|\nu|=\Vert \mu\circ R_k \Vert \leq |\mu| \Vert R_k \Vert = |\mu| < \sum_{i=1}^n |\mu(A_i)|+\varepsilon$ for
some sequence $(A_i)$ of disjoint closed subsets of $K$. Hence,
\begin{align*}
|\nu| &< \sum_{i=1}^n \left|\int_K 1_{A_i} \, d\mu \right|+\varepsilon = \sum_{i=1}^n |P_{A_i}(1)|+\varepsilon  \\
      &\leq \sum_{i=1}^n \Vert P_{A_i} \Vert +\varepsilon \\
      &\leq \Vert P \Vert +\varepsilon ,
\end{align*}
this last inequality by Lemma 1.2 of \cite{CaLaZa}.
\end{proof}

Note that in the proof we have shown that the natural and integral norms coincide on the
space of orthogonally additive $k$-homogenous polynomials on $C(K)$.
We have also shown that any orthogonally additive function in  $\mathcal{H}_b (C(K))$ is,
in fact, in  $\mathcal{H}_{bI} (C(K))$, which means that its restriction to $nB_{C(K)}$ belongs to
$\mathcal{H}_I (nB_{C(K)})$ for all $n$, (see \cite{DiGaMaZa}).

The converse of Proposition~\ref{bOA implies integral} does not hold even for homogeneous polynomials:
the condition for an integral
polynomial $P$ to be orthogonally additive is that it factor through $R_k$, i.e.
$$
\sum_i a_i x_i^k =0 \text{ on }K \Rightarrow \sum_i a_i P(x_i) =0.
$$
Now consider the analogous situation for an integral holomorphic map
$f\colon B^\circ_{C(K)} \rightarrow \mathbb{C}$. Take $X$ the closed subspace of
$C(B_{{C(K)}'})$ spanned by $\{\frac{\widehat{x}}{1 - \widehat{x}} : \Vert x \Vert < 1 \}$,
and $R\colon X \rightarrow C(K)$ given by
$R \left( \frac{\widehat{x}}{1 - \widehat{x}} \right) = \frac{x}{1 - x}$.
Then $R$ is again well-defined and continuous, and by composition with $R$, all maps of the
form
$$
f(x)=\int_K \frac{x}{1 - x} \, d\mu
$$
are integral (and $f(0)=0$). However, not all orthogonally additive bounded-type holomorphic
functions are of this form. To obtain such a representation, one
would have to construct a measure over $K$, and here the obstruction is algebraic: the spaces $R_k(X_k)$
are not `independent' in $C(K)$. The following example shows that for non-homogeneous polynomials
there is no such representation with any holomorphic $\Phi \colon \mathbb{C} \rightarrow \mathbb{C}$.

\begin{example}
Take $K$ to be the closed unit disc in $\mathbb{C}$. There is no entire function $\Phi \colon \mathbb{C}\rightarrow \mathbb{C}$
such that $h(x)=\Phi \circ x$ factors all degree two orthogonally additive polynomials over $C(K)$.
\end{example}

\begin{proof} Suppose there were such a function $\Phi$, i.e., given an orthogonally additive polynomial $P$,
there is a measure $\mu$ such that
$$
P(x)=\int_K \Phi \circ x \, d\mu \quad \quad \text{ for }x\in C(K).
$$
Write the Taylor series of $\Phi$ at 0:
$$
\Phi (z)=\sum_{k=0}^\infty c_k z^k
$$
and choose $a$ and $b$ in $K$ such that $c_2 a^2 \neq c_1 b^2$. Define $P\colon C(K)\rightarrow \mathbb{C}$
by $P(x)=x(a)+x(b)^2$. Clearly, $P$ is orthogonally additive, so let $\mu$ be as above.
Now for any $\lambda$,
$$
x(a)\lambda + x(b)^2\lambda^2=P(\lambda x)=\int_K \Phi \circ \lambda x \, d\mu =\sum_{k=0}^\infty c_k \int_K x^k \, d\mu \, \lambda^k
$$
and thus
$$
x(a)=c_1 \int_K x \, d\mu , \quad \text{ and } \quad x(b)^2=c_2\int_K x^2 \, d\mu \text{ for }x\in C(K).
$$
However, if we take $x(z)=\frac{c_2}{c_1}z^2$ and $y(z)=z$ (and thus $c_1 x=c_2 y^2$),
$$
\frac{c_2}{c_1}a^2=x(a)=c_1 \int_K x \, d\mu = c_2\int_K y^2 \, d\mu =y(b)^2=b^2,
$$
which contradicts our choice of $a$ and $b$.
\end{proof}

\section{Characterization of orthogonally additive holomorphic functions of bounded type}

For any Borel measure $\mu$ on K, we will write $f<<\mu$, if $\mu (A)=0$ implies $f_A=0$.
Also, we will say that $h\colon C(K) \rightarrow L^1 (\mu)$ is a power series function if
$$
h(x)=\sum_{k=1}^\infty g_k x^k ,
$$
with $g_k$'s in $L^1 (\mu)$. Note that this is a very special type of holomorphic function,
as its definition uses the algebra structure of $C(K)$.
We then have the following.

\begin{theorem}
Given $f\in \mathcal{H}_b(C(K))$, orthogonally additive, $f<<\mu$
if and only if there exists a power series function $h\in \mathcal{H}_b(C(K), L^1 (\mu))$
such that
$$
f(x)=\int_K h(x)(t) \, d\mu (t) \quad \quad \text{ for all }x\in C(K).
$$
\end{theorem}

\begin{proof}
Consider the Taylor series of $f$ about 0: $f=\sum_k P_k$. Note that $\overline{f}=\sum_k \overline{P_k}$
in $C(K)''$. Thus for any Borel set $A$,
$$
f_A(x)=\overline{f}(1_A x)=\sum_{k=0}^\infty \overline{P_k}(1_A x)=\sum_{k=0}^\infty P_{kA}(x).
$$
Since $\Vert P_{kA} \Vert \leq \Vert P_k \Vert$ by \cite[Lemma 1.2.]{CaLaZa}
$f_A \in \mathcal{H}_b(C(K))$

$\Rightarrow$) Now suppose $f<<\mu$. Since $f_A=0$ if and only if $P_{kA}=0$ for all $k$, we have
$P_k<<\mu$ for all $k$. Also, since $f$ is orthogonally additive, so is $P_k$ for all $k$.
Therefore, by \cite{BLLl} \cite{PGV}, there is a measure $\mu_k$ such that
$$
P_k(x)=\int_K x^k \, d\mu_k ,
$$
but since $P_k<<\mu$, $\mu_k<<\mu$: indeed, if $\mu(A)=0$,
$$
\mu_k(A)=\int_K 1_A \, d\mu_k=\overline{P_k}(1_A)=P_{kA}(1)=0,
$$
(the second equality by \cite[Corollary 2.1.]{CaLaZa}).
Now, by the Radon-Nikodym theorem there is a $g_k\in L^1(\mu)$ such that
$$
P_k(x)=\int_K x^k g_k \, d\mu.
$$
We prove that $\Vert g_k \Vert_{L^1(\mu)} \leq \Vert P_k \Vert$ for each $k$: choose a representative
$g_k$ and consider the function
$$
t \mapsto \left\{
\begin{aligned}
\left(\frac{\overline{g_k}(t)}{|g_k(t)|}\right)^{1/k}, & \text{ if } g_k(t)\neq 0, \\
0 \quad \quad \quad , & \text{ if } g_k(t)=0,
\end{aligned}\right.
$$
where, in taking $k$-th root we have chosen any branch. This function, though not continuous,
is Borel measurable, and has image in the closed unit disc. It is an element of $C(K)''$.
Using \cite[Corollary 2.1.]{CaLaZa} again we have
\begin{align*}
\Vert g_k \Vert_{L^1(\mu)}
& =\int_K |g_k|\, d\mu
=\int_K \left[\left(\frac{\overline{g_k}(t)}{|g_k(t)|}\right)^{1/k}\right]^k g_k \, d\mu \\
& =\overline{P_k}\left(\left(\frac{\overline{g_k}(t)}{|g_k(t)|}\right)^{1/k}\right)
\leq \Vert \overline{P_k} \Vert = \Vert P_k \Vert,
\end{align*}
this last equality by \cite{DG}.
Now define $h\colon C(K) \rightarrow L^1(\mu)$ by
$$
h(x)=\sum_{k=1}^\infty g_k x^k.
$$
This series converges absolutely:
\begin{align*}
\left\Vert \sum_{k=0}^\infty g_k x^k \right\Vert_{L^1(\mu)}
& \leq \sum_{k=0}^\infty \Vert g_k x^k \Vert_{L^1(\mu)} \\
& = \sum_{k=0}^\infty \int_K |g_k| |x|^k \, d\mu
\leq \sum_{k=0}^\infty \Vert g_k \Vert_{L^1(\mu)} \Vert x \Vert_\infty^k \\
& \leq \sum_{k=0}^\infty \Vert P_k \Vert \Vert x \Vert^k < \infty \quad \text{ for all } x.
\end{align*}
Thus $h$ is a power series function and is clearly bounded on bounded subsets of $C(K)$. So
$h\in \mathcal{H}_b(C(K), L^1 (\mu))$, and
$$
f(x)=\sum_{k=0}^\infty P_k (x)=\sum_{k=0}^\infty \int_K g_k  x^k \, d\mu= \int_K \sum_{k=0}^\infty g_k  x^k \, d\mu=
\int_K h(x) \, d\mu .
$$
$\Leftarrow$) If $\mu(A)=0$,
$$
f_A(x)=\overline{f}(1_A x)=\int_K \overline{h}(1_A x) \, d\mu = \int_K h(x) 1_A \, d\mu = \int_A h(x) \, d\mu =0,
$$
so $f<<\mu$. The second equality holds because $h$ is a power series function of bounded type.
Indeed, $ f=\sum_k P_k =\int_K h \, d\mu$. Since $h$ is a power series, $P_k (x)= \int_K g_kx^k  \, d\mu$.
The proof of \cite[Corollary 2.1]{CaLaZa} then shows that $\overline{P_k} (1_A x)=\int_A g_k x^k  \, d\mu$,
and the result follows.
\end{proof}

\begin{proposition}
If $f$ is orthogonally additive, then for some Borel measure $\mu$, $f<<\mu$.
\end{proposition}

\begin{proof}
Say $f=\sum_k P_k$ is the Taylor expansion of $f$ at zero.
Since $f$ is orthogonally additive, so are all the $P_k$'s. Thus there are
Borel measures $\mu_k$ for which
$$
P_k(x)=\int_K x^k \, d\mu_k .
$$
Now, $P_{kA}(x)=\overline{P_k}(1_A x)=\int_A x^k \, d\mu_k$, so $P_{kA}(1)=\mu_k(A)$ and $|\mu_k|\leq \Vert P_k \Vert$.
Therefore $(|\mu_k|)$ is summable, and we can define $\mu=\sum_k |\mu_k|$. For any $A$,
$\mu(A)=\sum_k |\mu_k|(A)$, so $\mu_k<<\mu$ for all $k$. This implies that
$P_k<<\mu$ for all $k$, and hence $f<<\mu$.
\end{proof}

We have the following characterization of orthogonally additive holomorphic functions of bounded type.

\begin{theorem}
$f\in \mathcal{H}_b(C(K))$ is orthogonally additive if and only if there is a Borel measure $\mu$ and a
power series function $h\in \mathcal{H}_b(C(K), L^1 (\mu))$ such that
$$
f(x)=\int_K h(x)(t) \, d\mu (t) \quad \quad \text{ for all }x.
$$
\end{theorem}

\begin{proof}
$\Rightarrow$) Immediate from the preceding results.

$\Leftarrow$) Note that all power series functions are orthogonally additive, for if $x$ and $y$ are orthogonal,
$(x+y)^k=x^k+y^k$ for all $k\geq 1$, so
$$
h(x+y)=\sum_{k=1}^\infty g_k (x+y)^k = \sum_{k=1}^\infty g_k (x^k+y^k) =
\sum_{k=1}^\infty g_k x^k + \sum_{k=1}^\infty g_k y^k = h(x)+h(y).
$$
Now if $f(x)=\int_K h(x) \, d\mu$, orthogonal additivity of $f$ follows from linearity of the integral.
Also, if $\Vert x \Vert \leq c$,
$$
|f(x)|\leq \sum_{k=1}^\infty \Vert g_k \Vert_{L^1 (\mu)} c^k,
$$
so $f$ is bounded on bounded sets.
\end{proof}

\end{document}